\documentclass[12pt]{article}
\usepackage[cp1251]{inputenc}
\usepackage[english]{babel}
\usepackage{amsmath}
\usepackage{amssymb}
\usepackage{amsfonts}
\usepackage{graphicx}
\usepackage{mathrsfs}
\newtheorem{lemma}{Lemma}
\newtheorem{theorem}{Theorem}
\newtheorem{corollary}{Corollary}
\numberwithin{equation}{section}
\numberwithin{theorem}{section}
\numberwithin{lemma}{section}

\begin{document}
\thispagestyle{empty}
\begin{center}
\bf \Large {Distortion and critical values \\of the finite Blaschke product}\end{center}\begin{center}
\bf V.N. Dubinin
\end{center}

\begin{quote}
\small \textbf {Abstract.}  We establish a sharp upper bound for the absolute value of the derivative of the finite Blaschke product, provided that the critical values of this product lie in a given disk.\end{quote}

\begin{quote}
\end{quote}
\section{Introduction}
\hspace{0.7cm}The inequalities for the absolute values of the derivative of a complex polynomial, taking into account its critical values, are of natural interest in the theory of multivalent functions.  A certain influence on the study of such inequalities was made by the well-known Erd\H{o}s conjecture about the maximum modulus of the derivative on a connected lemniscate [9-11].  The impact of the critical values of the polynomial on various type of distortions was considered in [4, 5].  In particular, the following result was obtained in [5]: if all critical values of the polynomial 
$$P(z)=c_0+c_1z+...+c_nz^n,\;\;c_n\neq 0,\;n\ge 2,$$
lie in the disk $|w|\le 1$ then
\begin{equation}|P'(z)|\le 2^{\frac{1-n}{n}}|c_n|^{\frac{1}{n}}T'_n(T_n^{-1}(|P(z)|))
\end{equation}
for every $z$. Equality holds in (1.1) for $P=T_n$ and all real $z,\;|z|\geqslant \cos(\pi/(2n)).$ Here $T_n(z)=2^{n-1}z^n+...$ is the Chebyshev polynomial of the first kind of degree $n,$ and $T^{-1}_n(\cdot)$ is the continuous branch of its inverse function defined on the ray $[0,+\infty]$ and taking this ray onto the ray $[\cos(\pi/(2n)),+\infty]$ (see also [8]). The purpose of this article is to establish an inequality similar to (1.1) for the finite Blaschke products
$$B(z)=\alpha\prod_{k=1}^n\frac{z-z_k}{1-\overline{z}_k z},\;\;z\in U:=\{z:\;|z|<1\},$$
$|\alpha|=1,\;|z_k|<1,\;k=1,...,n.$ 
These products and their applications have been studied by many authors (see, for example, [12, 14-16] and the references therein).  At the same time, problems associated with critical values have not been fully studied [16; 17, p. 365].  In a number of problems on finite Blaschke products, the extremal function is the so-called \\ 
----------------------------------------\\
\small{Vladimir Nikolaevich Dubinin, dubinin@iam.dvo.ru\\Far Eastern Federal University, st. Sukhanov, 8, 690950, Vladivostok, Russia\\
Institute for Applied Mathematics, FEBRAS, st. Radio, 7, 690041, Vladivostok, Russia\\
The work is supported by the Russian Foundation for Basic Research (Grant no. 20-01-00018)}  Chebyshev-Blaschke product [15].  We need a function $B_{n\tau}$ which, up to a linear fractional replacement of the argument, coincides with the Chebyshev-Blaschke product. For a fixed positive integer $n>1$ and a number $\varkappa,\;0<\varkappa<1,$ consider the rational function  
$Z\equiv Z_n(\zeta;\varkappa)$ defined parametrically as $$Z_n({\rm sn}(u;k);\varkappa):={\rm sn}\left(u\frac{{\bf K}(\varkappa)}{{\bf K}(k)};\varkappa\right),\;\;u\in \mathbb C;$$
the modulus $k$ is determined from the condition $${\bf K}'(k){\bf K}(\varkappa)=n{\bf K}'(\varkappa){\bf K}(k),\;\;0<k<1,$$ where
${\bf K}(\cdot)$ and ${\bf K}'(\cdot)$ are complete elliptic integrals of the first kind [1]. The function $Z,$ as well as the compositions of $Z$ with linear-fractional transformations both on the domain of the argument and on the range of $Z,$ are customarily called {\it{Zolotarev fractions}} [2, 3]. The role of Zolotarev fractions in the rational approximation theory and calculating electrical filters is well known [1]. The function $B_{n \tau}$ is defined as the composition
$$B_{n \tau}(z)=\Phi\left(Z\left(\frac{z-1}{z+1};\varkappa\right)\right),\;\;z\in U,\;1<\tau<\infty,$$
where
$$\Phi(v)=\frac{v\sqrt{\varkappa}+1}{1-v\sqrt{\varkappa}},\;\;\;\;\;\sqrt{\varkappa}=\frac{\sqrt{\tau}-1}{\sqrt{\tau}+1}.$$
All zeros of the function $B_{n\tau}$ are simple, and its zeros are located on the interval $[(k-1)/(k+1),0]$ of the real axis. Everywhere below, $\beta$ denote the largest zeros of the function $B_{n\tau}$. The description of the function $B_{n\tau}$ will be provided in Section 2. Within this setting we can already formulate our main result. 
\begin{theorem}
Let $f$ be a finite Blaschke product of degree $n\ge 2,$ and let all critical values of $f$ lie in the disk $|w|\leqslant \lambda,\;0<\lambda<1.$ Then for any point $z\in U$
\begin{equation}
(1-|z|^2)|f'(z)|\leqslant (1-|B^{-1}_{n\tau}(|f(z)|)|^2)|B'_{n \tau}(B^{-1}_{n\tau}(|f(z)|))|,
\end{equation}
where $\tau=\lambda^{-2}$ and $B^{-1}_{n\tau}(|f(z)|)$ is taken in the interval $[(1+\beta)/(1-\beta),1].$ Equality is attained in {\rm (1.2)} for the composite functions of the form $B_{n \tau}\circ\varphi,$ where $\varphi$ is an arbitrary conformal automorphism of the disk $U$ such that $\varphi(z)\in [(1+\beta)/(1-\beta),1].$
\end{theorem}

Thus, we obtain a sharp upper bound for the absolute value of the derivative $f'$ at any point $z$ of the disk $U.$ This bound depends on $|f(z)|.$ Theorem 1.1 is proved in Section 4 by using the symmetrization method [7], in which the result of symmetrization is located on the Riemann surface of the function inverse to a Chebyshev polynomial (see Section 3). In Section 5 we give some corollaries of the theorem 1.1 which are of independent interest. 

\section{Zolotarev fractions and Chebyshev polynomial}
\hspace{0.7cm}Throughout what follows, a Riemann surface is a surface  $\mathscr R$ glued from finitely or countably many domains in the extended complex plane so that the following conditions are satisfied: each point in  $\mathscr R$  projects onto a point in one of the glued domains; each point in  $\mathscr R$  has a neighbourhood which is a univalent disk or a multivalent disk with the unique ramification point at the centre of the disk (see [13], Pt. 3 for details). When there can be no misunderstanding, we will not distinguish the plane domains before gluing (which identifies some parts of the boundaries of these domains) and after it (when they become subdomains of $\mathscr R$). 

For a meromorphic function  $f$ denote by $\mathscr R (f)$ the Riemann surface of the function inverse to the function $f.$ The function $Z$ maps the sphere $\overline{\mathbb C}_\zeta$ onto Riemann surface $\mathscr R(Z),$ lying over the sphere $\overline{\mathbb C}_v.$ One of the representations of this surface is as follows. Let $G_1$ be the $v$-plane cut along the union of the rays $g^{-}:=[-\infty,-1/\varkappa]\bigcup [1/\varkappa,+\infty],$ and let $G_2,...,G_{n-1}$ be copies of the $v$-plane cut along $g^{-}$ and the interval $g^{+}:=[-1,1];$ finally, let $G_n$ be the $v$-plane cut along $g^{-}$ if $n$ is even and along $g^{+}$ if $n$ is odd. We obtain the Riemann surface $\mathscr R(Z)$ by gluing together the domains $G_k,\;k=1,...,n,$ in the following way. The domain $G_1$ is glued with $G_2$ cross-wise along the sides of the $g^{-}$-cuts. The domain $G_2$ is glued with $G_3$ along the sides of the $g^{+}$-cuts and so on. The domain $G_{n-1}$ is glued with $G_n$ cross-wise along the sides of the $g^{-}$-cuts if $n$ is even and along the $g^{+}$-cuts if $n$ is odd. Regarded as subsets of the surface $\mathscr R(Z),$ the domains $G_k$ involved in the procedure will be denoted by $\mathscr G_k,\;k=1,...,n,$ respectively. We treat the function $Z$ as a map of the sphere $\overline{\mathbb C}_\zeta$ on the Riemann surface $\mathscr R(Z)$ under which $Z([-1,1])\subset \mathscr G_1 (\rm{pr}\,{\it Z}([-1,1])=[-1,1]).$ This description of the Riemann surface $\mathscr R(Z)$ is easy to obtain by using properfies of the elliptic sine and the Riemann-Schwarz symmetry principle for conformal mappings. The compact Riemann surface $\mathscr R(Z)$ is a schlichtartig surface, and it has a finite number of sheets. Therefore, the function $Z$ is rational function of degree $n$ [13, Chapter 8, \S 10, Theorem 4]. We note that the function $Z$ maps the left half-plane onto the left half-plane. Taking  into account the properties of the linear-fractional mappings, we coclude that the function $B_{n\tau}$ is a finite Blaschke product, and the representation of the  Riemann surface $\mathscr R(B_{n\tau})$ is as follows. Let $H_1$ be the disk $U_w$ cut along the interval $h^{-}:=[-1,-1/\sqrt{\tau}],$   and let $H_2,...,H_{n-1}$ be copies of the disk $U_w$ cut along the intervals $h^{-}$ and $h^{+}:=[1/\sqrt{\tau},1];$ finally, let $H_n$ be the disk $U_w$ cut along $h^{-}$ if $n$ is even and along $h^{+}$ if $n$ is odd. We obtain the Riemann surface $\mathscr R(B_{n\tau})$ by gluing together the domains $H_k,\;k=1,...,n,$ in the following way. The domain $H_1$ is glued with $H_2$ cross-wise along the sides of the $h^{-}$-cuts. The domain $H_2$ is glued with $H_3$ along the sides of the $h^{+}$-cuts and so on. The domain $H_{n-1}$ is glued with $H_n$ cross-wise along the sides of the $h^{-}$-cuts if $n$ is even and along $h^{+}$-cuts if $n$ is odd. Regarded as subsets of the surface $\mathscr R(B_{n\tau}),$ the domains $H_k$ will be denoted by $\mathscr H_k,\;k=1,...,n,$ respectively. We treat the function $B_{n\tau}$ as a map of the disk $U_z$ on the Riemann surface $\mathscr R(B_{n\tau})$ under which $B_{n\tau}([(1+\beta)/(1-\beta),1))\subset \mathscr H_1({\rm{pr}}\,B_{n\tau}([(1+\beta)/(1-\beta),1))=[0,1)).$

The Riemann surface $\mathscr R(\sqrt{\tau}B_{n\tau})$ approaches the surface $\mathscr R(T_n)$ when $\tau\to \infty.$ We need a description of the surface  $\mathscr R(T_n)$ taken from [7]. Let $D_1$ be the $w$-plane cut along the ray $d^{-}:=[-\infty,-1],$ and let $D_2,...,D_{n-1}$ by copies of the $w$-plane cut along the rays $d^{-}$ and $d^{+}:=[1,+\infty];$ finally, let $D_n$ be the $w$-plane cut along the ray $d^{-}$ if $n$ is even or along $d^{+}$ if $n$ is odd. We obtain the Riemann surface $\mathscr R(T_n)$ by gluing together the domains $D_k,\;k=1,...,n,$ as follows. We glue $D_1$ with $D_2$ cross-wise along the sides of the $d^{-}$-cuts. The domain $D_2$ is glued with $D_3$ along the sides of the $d^{+}$-cuts and so on. The domain $D_{n-1}$ is glued with $D_n$ cross-wise along the sides of the $d^{-}$-cuts if $n$ is even and along the $d^{+}$-cuts if $n$ is odd. Regarded as subsets of the surface $\mathscr R(T_n),$ the domains $D_k$ will be denoted by $\mathscr D_k,\;k=1,..,n,$ respectively. Let $\mathscr L$ denote the ray lying on the sheet $\mathscr D_1$ over the ray $[0,+\infty].$ Then $T_n([\cos(\pi/(2n)),+\infty])=\mathscr L.$ We consider the surface $\mathscr R(\sqrt{\tau}B_{n\tau})$ as a subset of the surface $\mathscr R(T_n).$ The representation of the Riemann surface $\mathscr R(T_n)$ is used to determine the symmetrization in the next section.
\section{Symmetrization}
\hspace{0.7cm}This section provides information on the symmetrization [7] to the extent it is necessary to prove  Theorem 1.1. Let $\gamma(\rho)=\{w:\;|w|=\rho\},\;0\leqslant\rho\leqslant \infty.$ Then $\mathfrak R_n,\;n\geqslant 1,$ will denote the class of Riemann surfaces $\mathscr R$ over the complex $w$-sphere which satisfy the following conditions:

1) taking account of multiplicities, the total linear measure of any system of arcs on $\mathscr R$ which lies over an arbitrary circle $\gamma(\rho),\;0<\rho<\infty,$ has the estimate $2\pi n\rho;$

2) for $1\leqslant\rho<\infty,$ any closed Jordan curve on $\mathscr R$ over a circle $\gamma(\rho)$ which does not pass throught ramification points of $\mathscr R$ covers this circle with multiplicity $n.$

We now take an arbitrary surface $\mathscr R$ of class $\mathfrak R_n$ and proceed to the definition of the circular symmetrization $\rm{Sym}$ of sets and condensers on $\mathscr R$ [7]. Let ${\mathscr B}$ be an open set in ${\mathscr R}.$ Then symmetrization ${\rm Sym}$ transforms ${\mathscr B}$ into a subset ${\rm Sym}\hspace{0,5mm}{\mathscr B}$ of ${\mathscr R}(T_n)$ with the following properties. Fix some $\rho,\;0\leqslant \rho \leqslant \infty.$ If no points in ${\mathscr B}$ lie over the circle $\gamma(\rho),$ then no points in ${\rm Sym}\hspace{0,5mm}{\mathscr B}$ lie over it either. If ${\mathscr B}$ covers $\gamma(\rho),\;1\leqslant \rho\leqslant\infty,$ with multiplicity $n,$ then ${\rm Sym}\hspace{0,5mm}{\mathscr B}$ also covers $\gamma(\rho)$ with  multiplicity $n.$ If ${\mathscr B}$  covers $\gamma(\rho),\;0\leqslant \rho <1,$ with multiplicity $l\leqslant n,$ then the part of ${\rm Sym}\hspace{0,5mm}{\mathscr B}$ over $\gamma(\rho)$ consists of $l$ circles lying on the sheets ${\mathscr D}_1,...,{\mathscr D}_l.$ In the other cases, for $1\leqslant\rho<\infty$ the part of ${\rm Sym}\hspace{0,5mm}{\mathscr B}$ lying over $\gamma(\rho)$ is an open arc\footnote[1]{For $\rho>1,$ this is an open Jordan arc, but for $\rho=1$ it can contain self-tangency points.} on ${\mathscr R}(T_n)$ with midpoint on the ray ${\mathscr L}$ and with linear measure equal to the measure of ${\mathscr B}(\rho):=\{W\in{\mathscr B}:\;|{\rm pr}W|=\rho\}.$ For $0<\rho<1$ the part of ${\rm Sym}\hspace{0,5mm}{\mathscr B}$ over $\gamma(\rho)$ is a union of $m$ circles $\Gamma_1,...,\Gamma_m,\;0\leqslant m\leqslant n-1,$ and an open arc $\Gamma_{m+1}$ such that $\Gamma_k=\Gamma_k({\mathscr B},\rho)\subset {\mathscr D}_k,\;k=1,...,m+1;$ the total linear measure of these curves is equal to the measure of ${\mathscr B}(\rho),$ and the midpoint of $\Gamma_{m+1}$ lies over $(-1)^m\rho.$ Here the number $m$ of the circles depends on the measure of ${\mathscr B}(\rho).$ If this measure is less that $2\pi\rho,$ then necessarily $m=0,$ and there are no full circles.

The result ${\rm Sym}\hspace{0,5mm}{\mathscr E}$ of the symmetrization of a closed set ${\mathscr E}\subset {\mathscr R}$ also lies on ${\mathscr R}(T_p)$ and is defined as follows. Fix some $\rho,\;0\leqslant \rho\leqslant \infty.$ If no points in the set ${\mathscr E}$ lie over $\gamma(\rho),$ then ${\rm Sym}\hspace{0,5mm}{\mathscr E}$ contains no points over this circle either. If ${\mathscr E}$ covers $\gamma(\rho),\;1\leqslant\rho\leqslant\infty,$ with multiplicity $n,$ then ${\rm Sym}\hspace{0,5mm}{\mathscr E}$  also covers $\gamma(\rho)$ with multiplicity $n.$ If ${\mathscr E}$ covers $\gamma(\rho),\;0\le\rho< 1,$ with multiplicity $l\le n,$ then the part of ${\rm Sym}\hspace{0,5mm}{\mathscr E}$ over $\gamma(\rho)$ consists of $l$ circles in the sheets ${\mathscr D}_1,...,{\mathscr D}_l.$ Otherwise the part of ${\rm Sym}\hspace{0,5mm}{\mathscr E}$ lying over $\gamma(\rho),\;1\leqslant\rho< \infty,$ is a closed arc segment (that is, an arc with its endpoints) on ${\mathscr R}(T_n),$ with midpoint on the ray ${\mathscr L}$ and with linear measure equal to the measure of ${\mathscr E}(\rho):= \{W\in{\mathscr E}:\;|{\rm pr}W|=\rho\}$ (if the latter is equal to zero, then the corresponding arc segment is a point on ${\mathscr L}$). The part of ${\rm Sym}\hspace{0,5mm}{\mathscr E}$ over $\gamma(\rho),\;0<\rho< 1,$ is the union of $m$ circles $\Gamma_1,...,\Gamma_m,\;0\le m\le n-1,$ and a closed arc segment $\Gamma_{m+1}$ such that $\Gamma_k\subset{\mathscr D}_k,\;k=1,..,m+1,$ the total linear measure of these curves is equal to that of ${\mathscr E}(\rho)$ and the midpoint of $\Gamma_{m+1}$ lies over  
$(-1)^m\rho$ (if the measure in question is $2\pi\rho m,$ where $m$ is a nonnegative integer, then $\Gamma_{m+1}$ reduce to a point).

A {\it condenser on the surface ${\mathscr R}$} is an ordered pair of sets ${\mathscr C}=({\mathscr B},{\mathscr E}),$ where ${\mathscr B}$ is an open subset of ${\mathscr R}$ and ${\mathscr E}$ is a compact subset of ${\mathscr B}.$ We call ${\mathscr B}\setminus {\mathscr E}$ the {\it field of the condenser ${\mathscr C}$.} Now we set 
$${\rm Sym}\hspace{0,5mm}{\mathscr C}=({\rm Sym}\hspace{0,5mm}{\mathscr B},\;{\rm Sym}\hspace{0,5mm}{\mathscr E}).$$
The {\it capacity} ${\rm cap}\hspace{0,5mm}{\mathscr C}$ of the condenser ${\mathscr C}=({\mathscr B},{\mathscr E})$ is defined by
$${\rm cap}{\mathscr C}=\inf \int_{{\mathscr B}}|\bigtriangledown{\mathscr V}|^2 d\sigma,$$ where the infimum is taken over all {\it admissible} functions ${\mathscr V}:$ real-valued functions ${\mathscr V}$ which have compact support in ${\mathscr B},$ are equal to 1 on ${\mathscr E}$ and are locally Lipschitz in ${\mathscr B}.$ If there exists a function ${\mathscr P}$ which is continuous in $\overline{{\mathscr B}}$, equal to zero on $\partial {\mathscr B},$ to 1 on  ${\mathscr E},$ and harmonic in the field ${\mathscr B}\setminus{\mathscr E},$ then it is called the {\it potential function} of the condenser ${\mathscr C}$. Then by Dirichlet's principle $${\rm cap}\hspace{0,5mm}{\mathscr C}=\int_{{\mathscr B}\setminus{\mathscr E}}|\bigtriangledown{\mathscr P}|^2 d\sigma.$$  

The following is a special case of the central result of [7].
\vskip 2mm
\hspace{-0.4cm}{\bf Lemma 3.1.} {\rm (see [7, Theorem 1.1]).} {\it For each condenser ${\mathscr C}$ on a surface ${\mathscr R}$ 
\begin{equation}
{\rm cap}\hspace{0,5mm} {\mathscr C}\ge{\rm cap\hspace{0,5mm} Sym}\hspace{0,5mm} {\mathscr C}.
\end{equation}}

For condencers $\mathscr C$ with a connected field having a potential function, all cases of equality in (3.1) are known [7]. 
 
 \section{Proof of Theorem 1.1}
\hspace{0.7cm} Fix any point $z_0\in U$ such that $f'(z_0)\neq 0,$ and look at the condenser $$C(r)=(U,\{z:|z-z_0|\leqslant r\})$$
where $r>0$ is sufficiently small. Denote by $\mathscr C(r)$ the image of $C(r)$ under the map $F:=f/\lambda,$ which we view as lying on the Riemann surface $\mathscr R(F).$ The surface $\mathscr R(F)$ belongs to the class $\mathfrak R_n$ from Section 3. Indeed, the condition 1) is obvious. To verify condition 2), we consider an arbitrary fixed $\rho,\;1<\rho<1/\lambda,$ and  
 the sets $W_\rho=\{w:\;\rho<|w|<1/\lambda\},$ $V_\rho=F^{-1}(W_\rho).$ Let $F(V_\rho)$ be the image of $V_\rho$ on the surface $\mathscr R(F)$. Under the map $F,$ the points of the unit circle $|z|=1$ and only they transform to the points of the circle $|w|=1/\lambda.$ Therefore, there is only one boundary component of the set $F(V_\rho)$ lying over $|w|=1/\lambda$ and covering this circle $n$-fold. It follows in particular that the set $V_\rho$ is a domain. According to the hypothesis of Theorem 1.1, the surface $F(V_{\rho})$ $n$-fold covers the ring $W_{\rho}$ and does not contain ramification points over $W_{\rho}.$ Using the Hurwitz formula, we conclude that the domain $V_{\rho}$ has only two boundary components. Therefore, the circle $|w|=\rho$ is $n$-fold covered by a connected set, and thus condition 2) is satisfied. From the above it follows that the symmetrization Sym from Section 3 is applicable to the condenser $\mathscr C(r).$ Taking into account the conformal invariance of condenser capacity and Lemma 3.1, we see that
 \begin{equation}
{\rm cap}\hspace{0,5mm} C(r)={\rm cap}\hspace{0,5mm} {\mathscr C}(r)\geqslant{\rm cap\hspace{0,5mm} Sym}\hspace{0,5mm} {\mathscr C}(r).
\end{equation}
The condenser ${\rm Sym}\hspace{0,5mm} {\mathscr C}(r)$ has the form $(\mathscr B,\mathscr E(r)),$ where $\mathscr B$ is a part of the surface $\mathscr R(T_n)$ lying over the disk $|w|<1/\lambda,$ and a set $\mathscr E(r)$ is almost a disk on the sheet $\mathscr D_1\subset \mathscr R(T_n)$ of radius   
\begin{equation}
(r|f'(z_0)|/\lambda)(1+o(1)),\;\;r\to 0,
\end{equation} 
centered at a point $\omega$ on the ray $\mathscr L,$ ${\rm pr}\,\omega=|f(z_0)|/\lambda.$ In other words, a close set $\mathscr E(r)$ contains and is contained in a disk centered at the point $\omega$ and a radius in the form of (4.2). Let $\mathscr C_1(r)$ be the result of extension the condenser $\rm{Sym}\, {\mathscr C}({\it r})$ by a factor of $\lambda.$ The condenser $\mathscr C_1(r)$ is located on the surface $\mathscr R(B_{n\tau}),$ where $\tau=\lambda^{-2}.$ Finally, we denote by $C_1(r)$ the image of the condenser $\mathscr C_1(r)$ under the map $B^{-1}_{n\tau}.$ This condenser has the form
$$C_1(r)=(U,E(r)),$$
where $E(r)$ is almost a disk of radius 
$$\frac{r|f'(z_0)|(1+o(1))}{|B'_{n\tau}(B^{-1}_{n\tau}(|f(z_0)|))|},\;\;r\to 0,$$
centered at the point $B^{-1}_{n\tau}(|f(z_0)|).$ Here $B^{-1}_{n\tau}$ is the continuous branch of the fuction inverse to $B_{n\tau}$ mapping the half-open interval $[0,1)$ on the set $[(1+\beta)/(1-\beta),1).$ Since the capacity is conformally invariant, we have 
$${\rm cap\hspace{0,5mm} Sym}\hspace{0,5mm} {\mathscr C}(r)={\rm cap}\hspace{0,5mm} {\mathscr C}_1(r)={\rm cap}\hspace{0,5mm} C_1(r).$$

With (4.1) in view, we conclude that
\begin{equation}
{\rm cap}\hspace{0,5mm} C(r)\geqslant{\rm cap}\,C_1(r).
\end{equation}

The asymptotic formula for condenser capacities from (4.3), as $r\to 0,$ is well known (see, for example, the more general case [6, Theorem 2.1]). Applying this formula, we obtain from (4.3) the inequality (1.2), where $z=z_0.$ In the case when $f'(z_0)=0,$ the inequality (1.2) is trivial.

Now suppose that for a given point $z_0\in U$ the function $f$ is equal to $B_{n\tau}\circ\varphi,$ where $\varphi$ is a comformal automorphism of the disk $U$ such that $\varphi(z_0)\in [(1+\beta)/(1-\beta),1].$ For the critical points $z$ of the function $f,$ the absolute values of $f(z)$ are equal to $1/\sqrt{\tau}.$ In the above notation, we set $\lambda=1/\sqrt{\tau}.$ Let $\tilde{\mathscr E}(r)$ be a closed disk on the sheet $\mathscr D_1$ contered at the point $\omega$ with radius (4.2) and containing the set $\mathscr E(r).$ The symmetrization Sym does not change the condencer $(\mathscr B,\tilde{\mathscr E}(r)):$
$${\rm Sym}\hspace{0,5mm} (\mathscr B,\tilde{\mathscr E}(r))=(\mathscr B,\tilde{\mathscr E}(r)).$$
Using the monotonicity of the capacity we obtain
$${\rm cap}\hspace{0,5mm} C(r)={\rm cap}\hspace{0,5mm} {\mathscr C}(r)\leqslant{\rm cap}\hspace{0,5mm}(\mathscr B,\tilde{\mathscr E}(r))={\rm cap\hspace{0,5mm} Sym}\hspace{0,5mm} (\mathscr B,\tilde{\mathscr E}(r)).$$
Repeating the previous proof with the replacement of the condenser ${\rm Sym}\hspace{0,5mm} \mathscr C(r)$ by ${\rm Sym}\hspace{0,5mm} (\mathscr B,\tilde{\mathscr E}(r))$, we arrive at the inequality opposite to (4.3). This implies the inequality opposite to (1.2), which, taking into account what was proved earlier, gives equality in (1.2). This completes the proof of Theorem 1.1.
 \section{Corollaries}
 
\hspace{0.7cm} Let $f$ be a meromorphic function in the disk $U$ different from a constant and let $0\leqslant t\leqslant \infty.$ The set 
 $$L_f(t)=\{z\in U:\;|f(z)|=t\}$$
 is called a {\it lemniscate} of the function $f.$ The following statement is a solution of an analogue of the Erd\H{o}s problem for the complex polynomials (see [9], [10]).
\begin{corollary}
Let $f$ be a finite Blaschke product of degree $n\geqslant 2,$ such that the lemniscate $L_f(\lambda)$ is connected for some $\lambda,\;0<\lambda<1.$ Then 
$$\max\{(1-|z|^2)|f'(z)|:\;z\in L_f(\lambda)\}\leqslant |B'_{n\tau}(0)|,$$
where $\tau=\lambda^{-2}.$ Equality holds for $f=B_{n\tau},\;\tau=\lambda^{-2}.$
\end{corollary}
{\bf Proof.}
Applying the Hurwitz formula to the $n$-fold covering by $\mathscr R(f)$ of the ring $\lambda<|w|<1,$ we conclude that there are no ramification points of the surface $\mathscr R(f)$ lying over this ring. By Theorem 1.1, for any point $z$ on the lemniscate $L_f(\lambda)$ we have 
$$(1-|z|^2)|f'(z)|\leqslant(1-|B^{-1}_{n\tau}(\lambda)|^2)|B'_{n\tau}(B^{-1}_{n\tau}(\lambda))|=|B'_{n\tau}(0)|,$$
where $\tau=\lambda^{-2}.$ Note that the point $z=0$ belongs to the lemniscate $L_f(\lambda)$ when $f=B_{n\tau}.$ The corollary 5.1 is proved.
\begin{corollary}
Let $f$ be a finite Blaschke product of degree $n\geqslant 2,$ and let $f(0)=0.$ Suppose that all critical values of $f$ lie in the disk $|w|\leqslant \lambda,\;0<\lambda<1.$ Then
$$|f'(0)|\leqslant |f'_{n\tau}(0)|,$$ 
where $f_{n\tau}=B_{n\tau}\circ \varphi,\;\tau=\lambda^{-2}$ and 
$$\varphi(z)=\frac{(1-\beta)z+(1+\beta)}{(1-\beta)+(1+\beta)z}.$$
\end{corollary} 
{\bf Proof.}
Using Theorem 1.1 we obtain 
$$|f'(0)|\le \left[1-\left(\frac{1+\beta}{1-\beta}\right)^2\right]\left|B'_{n\tau}\left(\frac{1+\beta}{1-\beta}\right)\right|=|f'_{n\tau}(0)|.$$
This proves Corollary 5.2

To verify the next corollary of Theorem 1.1, we need the following property of the extremal function $f_{n\tau.}$
\begin{lemma}
The absolute value $|f'_{n\tau}(0)|$ is a continuous strictly decreasing function of the variable $\tau\in [1,+\infty)$ mapping the ray $[1,+\infty)$ on the half-interval $(0,1].$
\end{lemma}
{\bf Proof.}
The continuity follows from the analytic representation of the function $f_{n\tau.}$ Unfortunately, the definition of this function does not make it possible to prove the monotonicity of $|f'_{n\tau}(0)|.$ Note that the Riemann surface $\mathscr R(f_{n\tau})$ coincides with the surface $\mathscr R(B_{n\tau}).$ Consequently
$$|f'_{{n\tau}}(0)|=r(U,0)|f'_{n\tau}(0)|=r(\mathscr R(f_{n\tau}),f_{n\tau}(0))=r(\mathscr R(B_{n\tau}),W_0),$$
where the point $W_0:=f_{n\tau}(0)=B_{n\tau}\left(\frac{1+\beta}{1-\beta}\right)$ lies on the sheet $\mathscr H_1,\;{\rm pr}\,W_0=0,$
and $r(D,W)$ means the inner radius of a domain $D$ with respect to a point $W$ [6, Ch. 2.1]. Thus, it is necessary to prove a strict decrease of the quantity $r(\mathscr R(B_{n\tau}),W_0)$ on the ray $[1,+\infty).$ Let $1<\tau_1<\tau_2<\infty,$ and let $g(W)$ be a Green function of the domain $\mathscr R(B_{n\tau_2})$ with pole at the point $W_0\in\mathscr H_1(\tau_2),\;{\rm pr}\,W_0=0.$ Here $\mathscr H_1(\tau_2)$ is the sheet $\mathscr H_1$ on the surface $\mathscr R(B_{n\tau})$ when $\tau=\tau_2$ (see Section 2). We fix a number $t>0$ such that the set $\mathscr E:=\{W:\;g(W)\geqslant t\}$ belongs to the sheet $\mathscr H_1(\tau_2).$ Let $\mathscr X$ be the union of all possible cuts on the surface $\mathfrak R(B_{n\tau_2})$ lying over the segments $[-1/\sqrt{\tau_1},-1/\sqrt{\tau_2}],$ $[1/\sqrt{\tau_2},1/\sqrt{\tau_1}].$ We consider the set $\mathfrak R(B_{n\tau_2})\setminus \mathscr X$ as a subset of the surface $\mathscr R(B_{n\tau_1}).$ Denote by $\mathscr P$ the potential function of the condenser $\mathscr C_2:=(\mathscr R(B_{n\tau_2}),\mathscr E),$ and let $\mathscr V$ be a the restriction of $\mathscr P$ on $\overline{\mathscr R(B_{n\tau_2})}\setminus\mathscr X$ extended on $\overline{\mathscr R(B_{n\tau_1})}$ by continuity. Such extension is possible in according to the symmetry of the function $\mathscr P$ with respect to the real axis on the appropriate sheets. By Dirichlet's principle
 \begin{equation}
 {\rm cap}\hspace{0,5mm}{\mathscr C}_2=\int_{\mathscr R(B_{n\tau_2})}|\bigtriangledown{\mathscr P}|^2 d\sigma
=\int_{\mathscr R(B_{n\tau_1})}|\bigtriangledown{\mathscr V}|^2 d\sigma>{\rm Cap}\,\mathscr C_1,
\end{equation}
where $\mathscr C_1=(\mathscr R(B_{n\tau_1}),\mathscr E).$ Further, we use Theorem 2.7 from [6] which is easily carried over to the case of Riemann surfaces. The application of this theorem twice together with inequality (5.1) leads to the following chain of relations 
$$\frac{1}{2\pi}\log r(\mathscr R(B_{n\tau_2}),W_0)=\frac{1}{2\pi}\log r(\mathscr E\setminus(\partial\mathscr E),W_0)+\frac{1}{{\rm cap}\,\mathscr C_2}<$$ 
$$<\frac{1}{2\pi}\log r(\mathscr E\setminus(\partial\mathscr E),W_0)+\frac{1}{{\rm cap}\,\mathscr C_1}\leqslant \frac{1}{2\pi}\log r(\mathscr R(B_{n\tau_1}),W_0).$$
The proof is complete.

In [16], a lower bound for the quantity
$$\max\left\{\left|\frac{f(\zeta)}{\zeta f'(0)}\right|:\;f'(\zeta)=0\right\}$$
was found for all finite Blaschke products of degree $n\geqslant2$ such that $f(0)=0,\;f'(0)\neq 0.$ The following sharp lower bound complements the research [16] on the critical values of Blaschke products.
\begin{corollary} 
Let $f$ be a finite Blaschke product of degree $n\geqslant 2,$ and let $f(0)=0,\;f'(0)\neq 0.$ Then
 \begin{equation}
 \max\{|f(\zeta)|:\;f'(\zeta)=0\}\geqslant \frac{1}{\sqrt{\tau}},
 \end{equation}
 where $\tau,\;1<\tau<\infty,$ is the unique root of the equation 
 $$|f'_{n\tau}(0)|=|f'(0)|.$$
 Equality holds in {\rm (5.2)} for the function $f_{n\tau}$ for any $\tau,\;1<\tau<\infty.$ 
\end{corollary}  
{\bf Proof.}
 Note that $f_{n\tau}(0)=0,\;f'_{n\tau}(0)\neq 0,$ and for any $\tau,\;1<\tau<\infty,$ the critical values of the function $f_{n\tau}$ are equal to $\pm 1/\sqrt{\tau}$. Denote by $\lambda$ the left side of (5.2). Using Corollary 5.2 we obtain
 $$|f'(0)|\leqslant |f'_{n\tilde{\tau}}(0)|,$$
 where $\tilde{\tau}=\lambda^{-2}.$ By Lemma 5.1, there exists the unique number $\tau$ such that 
 $$|f'_{n\tau}(0)|=|f'(0)|\leqslant |f'_{n\tilde{\tau}}(0)|.$$
 Again, by Lemma 5.1
 $$\tau\geqslant \tilde{\tau}.$$
 Here $\lambda\geqslant 1/\sqrt{\tau}.$ The corollary is proved. 


\end{document}